\documentclass[a4paper,12pt
]{article}
\usepackage{amsmath}
\usepackage{amssymb}
\usepackage[dvipdfmx]{graphicx}
\usepackage{graphicx}
\usepackage{here}
\usepackage{mathpazo}

\newtheorem{thm}{Theorem}[section]

\newtheorem{prop}[thm]{Propoition}

\newtheorem{cor}[thm]{Corollary}
\newtheorem{rem}[thm]{Remark}
\newtheorem{df}[thm]{Definition}

\newtheorem{constr}[thm]{Construction}
\newcommand{\QED}[1]{{\hspace*{\fill}#1}$\square$}

\newcommand{\R}{{\mathbb R}}
\newcommand{\N}{{\mathbb N}}
\newcommand{\C}{{\mathbb C}}
\newcommand{\Z}{{\mathbb Z}}

\newcommand{\F}{{\mathcal F}}
\newcommand{\G}{{\mathcal G}}

\newcommand{\IE}{{\it i.e.}, }
\newcommand{\EG}{{\it e.g.}, }

\newcommand{\ds}{\displaystyle}
\newcommand{\ol}[1]{{\overline{#1}}}
\begin{document}
\begin{center}{
{\bf \Large Modification of Convex Ends to Cylindrical 
\\
and 
Symplectic Foliations}
\vspace{15pt}
\\ 
Yoshihiko MITSUMATSU\footnote{
Supported 
by Grant-in-Aid for 
Scientific Research (B) 22340015, 17H02845, 21H00985.
  
{\it Date}: February 5, 2022.   
 {\copyright}  2022 Y. Mitsumatsu.  

}
}
\end{center}
\begin{abstract}
We show that the natural symplectic structure on 
the Milnor fiber of an isolated singularity 
in complex three variables whose link fibers over the circle 
can be modified into one which is cylindrical at the end. 
As a consequence we see that the foliation of codimension one on $S^5$ 
which is adapted to the Milnor open book of $S^5$ associated 
with such a singularity 
admits a leafwise symplectic structure. 
The modification also enables us to construct 
certain closed symplectic 4-manifolds. 
\end{abstract}
\setcounter{section}{-1}

\tableofcontents
\section{\large Introduction}
In the previous work \cite{Mi2} the author constructed 
leafwise symplectic structures on foliations 
adapted to the Milnor open books for 
the simple elliptic singularities 
$\tilde{E}_6, \tilde{E}_7$, and $\tilde{E}_8$.  
The foliation adapted to $\tilde{E}_6$ is so called 
{\em Lawson's foliation} (see \cite{L}) which is known 
to be the first 
foliation of codimension one constructed on the five sphere. 
In this paper we refine the construction mainly 
from the point of view of open book decomposition 
and give a sufficient condition not only for Milnor open book 
but also for general open books supporting contact structures 
so that the adapted foliation exists and admits 
a leafwise symplectic structure.  
As  a consequence, 
for isolated hypersurface singularities of three complex variables, 
Lawson type (adapted) foliations of the Milnor open books 
of cusp singularities also admit leafwise symplectic structures, 
as well as simple elliptic ones on which we have already worked.   

The notion of leafwise symplectic foliation or 
foliation with leafwise symplectic structure is 
essentially equivalent to that of regular Poisson structure 
and is defined as follows.  
\begin{df}\label{df:leafwise symplectic structure}
{\rm (Leafwise symplectic structure. See \cite{Mi2})  $\,$ 
For a smoothly foliated manifold $(M,\F)$, a leafwise symplectic structure 
is a smooth 2-form $\omega$ on $M$ which restricts to 
a symplectic form of each leaf of $\mathcal F$. 

Often it is regarded as a smooth section of 2nd exterior power 
of the cotangent bundle to the foliation, \IE an element of 
$\Gamma^{\infty}(\bigwedge^{2}T^{*}\F)$ 
whose restriction to each leaf is symplectic on the leaf. 
}
\end{df}
\begin{rem}{\rm $\,$ 
The one in the sense of the first definition restricts to one of the second definition, 
the one in the sense of second definition always admits a lift to 
one in the sense of the first definition.  
Therefore in this article we do not distinct the two and shift from 
one to the other without mention. 

The first definition has an importance if we would like to impose 
on the leafwise symplectic form $\omega$ to be closed as a 
2-form on the ambient manifold $M$.  For a foliation of codimension one 
with a globally closed leafwise symplectic form, 
Donaldson's theory of approximately holomorphic geometry 
is applicable and we find very fine structures (see Mart\'inez-Torres  
\cite{Mt}).  On the spheres a compact leaf is an apparent obstruction.  

The second definition is completely equivalent to the notion 
of regular Poisson structures, \IE  
Poisson structures all of whose symplectic leave has the same dimension.  
From this point of view our work is paraphrased as 
constructions of regular Poisson structures.   
}
\end{rem}

In Section \ref{sec:Milnor} we review the basic materials 
concerning singularities, including the Milnor fibrations 
and open books of isolated singularities, 
cusp and simple elliptic singularities. 
The associated foliations (often called Lawson type foliations)
with open books and their relations with the contact structures 
are reviewed in Section \ref{sec:Lawson}.

The key of the construction 
of leafwise symplectic structure on a Lawson type foliation 
is the modification of symplectic structure 
of a Milnor fiber 
into one with  cylindrical end. 
The natural symplectic structure of Milnor fiber which is inherited 
from the Euclidean space has a convex end. 
In order to construct a leafwise symplectic structure 
on an associated foliation, the fibers should spiral around 
a compact symplectic leaf, and accordingly the fibers 
is required to admit a symplectically periodic end. 
For Milnor fibers of simple elliptic singularities and 
cusp singularities we show that the modification is 
possible in a rather simple way (Section \ref{sec:MainTheorem}, 
Theorem \ref{thm:MainTheorem}). 

The method of modifying symplectic structures 
in Section \ref{sec:MainTheorem}
yields some further applications. 
We give some constructions of  symplectic structures
in Section \ref{sec:ClosedSymplectic} 
or 
of $b$-symplectic structures in Section \ref{sec:other}  
on certain closed 4-manifolds. 
A closed symplectic 4-manifold is obtained by gluing 
two copies of the Milnor fibers of certain cusp singularities 
(Construction \ref{constr:K3}), 
while $b$-symplectic structures are obtained on 
the doubles of the Milnor fiber 
of simple elliptic or cusp singularities.   
This construction seems to have some relation 
with recent study of homological mirror symmetry 
due to Hacking and Keating (\cite{K}, \cite{HK}).  
A closer look at these closed symplectic 4-manifolds 
is passed to a forthcoming paper \cite{KKMM}. 
\vspace{5pt}

\noindent
{\bf Acknowledgement}\quad The author is grateful to Atsuhide Mori 
for his insights and discussions as well as to 
Naohiko Kasuya for detailed discussions. He also is grateful 
to l'Unit\'e de Math\'ematique de l'Ecole Normale Sup\'erieure de Lyon 
because the main part of this work was done 
during the author's sabbatical stay in Lyon in 2012.

\section{Milnor open book 
and associated symplectic structure on pages}\label{sec:Milnor}
\subsection{Milnor fibration and open book decomposition}

First we recall the construction 
of exact symplectic open book 
decomposition from a Milnor fibration of an isolated singularity. 
Let $f(z_0, ... , z_n)$ be a polynomial or holomorphic function 
in $n+1$ complex variables 
and assume that the origin is an isolated critical point $\{f=0\}$.

Let $S^{2n+1}_{\rho_{0}}$ denote the sphere in $\C^{n+1}$ 
of radius $\rho_{0}>0$ with center at the origin.  
The most classical {\it Milnor fibration} 
is defined 
as 
$S^{2n+1}_{\rho_{0}} \setminus f^{-1}(0) \to S^1$ by taking 
$\mathrm{arg} f(z) =f(z)/\vert f(z) \vert$ for 
$z\in S^{2n+1}_r \setminus f^{-1}(0)$. 
Then for a small enough $\varepsilon_{0}>0$, 
the restriction of this fibration 
to a compact manifold 
$M=S^{2n+1}_{\rho_{0}} \cap \{ z=(z_0, \dots z_n)\,;\,
\vert f(z)\vert\geq \varepsilon_{0}\}$ 
with boundary 
is also called the Milnor fibration of the singularity. 
$U=S^{2n+1}_r\setminus M$ is a tubular neighborhood of the 
link $N=f^{-1}(0)\cap S^{2n+1}_r$ in the sphere $S^{2n+1}_r$. 
As the normal bundle of the link $N$ is trivialized by the value 
of $f$, $U$ is diffeomorphic to the product 
$N\times \mathrm{int} D^2_{\varepsilon_{0}}$ where 
$D^2_{\varepsilon_{0}}=\{w\in\C\,;\,\vert w \vert \leq \varepsilon_{0}\}$. 
The common boundary $W=\partial \overline{U} = \partial M$ 
is diffeomorphic to $S^1 \times N$. 
The decomposition $S^{2n+1}_r = \overline{U}\cup M$ 
together with the fibration $\varphi=f/\vert f \vert : M \to S^1$ 
is the Milnor open book decomposition of the sphere $S^{2n+1}$ 
associated with the isolated singularity. 
The binding is nothing but the link $N$. 

If we choose $\rho_{0}>0$ 
small enough carefully, the real-hypersurface  
$\{\mathrm{arg} f = \theta\}$ for each $\theta\in S^1$ is transverse to 
$S^{2n+1}_{\rho'}\setminus f^{-1}(0)$ for $0<\rho'\leq \rho_{0}$ so that 
the gradient flow of 
$\rho^2=\sum_{j=0}^n \vert z_j\vert^2$ 
on the real-hypersurfce
$\{\mathrm{arg} f=\theta\}$ 
gives rise to a diffeomprphism between 
$F_w=\{f(z)=w \} \cap D^{2(n+1)}_{\rho_{0}}$ and 
$\varphi^{-1}(\theta) \subset M \subset S^{2n+1}_{\rho_{0}}$ 
for $w=\varepsilon_{0}e^{i\theta}$.  
Therefore the fibration 
$\varphi : M \to S^1$ can be identified with  
$\ds 
\cup_{\vert w \vert=\varepsilon_{0}}F_w \to S^1$, 
which is also called a Milnor fibration.  
In this article, we mainly use this formulation of Milnor fibration 
and thus it is also denoted by $\varphi : M \to S^{1}$.  
When the absolute value $|w|$ is well-understood or 
a Milnor fiber $F_{w}$ is regarded as a page of the open book, 
$F_{w}$ is also denoted by $F_{\theta}$ where $\theta=\mathrm{arg}\,w$.

\subsection{Exact symplectic structures on Milnor pages 
and open book decompositions}
The canonical symplectic structure 
$\ds \omega_0=\sum_{j=0}^n dx_j\wedge dy_j$ on $\C^{n+1}$ 
with its canonical primitive (or a Liouville form) 
$\ds \lambda_0=\frac{1}{2}\sum_{j=0}^n 
(x_jdy_j - y_jdx_j)$ on $\C^{n+1}$ 
restricts to each Milnor fiber $F_w$  to be 
an exact symplectic structure. 
We slightly modify this structure by isotopying 
the embedding  $F_w \hookrightarrow \C^{n+1}$'s.

Let us identify the compact region 
$CM=\{\vert f \vert \leq \varepsilon_0\} \cap 
(D^{2(n+1)}_{\rho_{0}} \setminus \mathrm{int} D^{2(n+1)}_{{\rho_{0}}/2}) 
\subset \C^{n+1}$
with the product $[{\rho_{0}}/2, {\rho_{0}}]
\times D^2_{\varepsilon_0}\times N$ as follows.  
The first coordinate $\rho$ is the distance from the origin, 
the second is the value of $f$ in 
$D^2_{\varepsilon_0}\subset\C$. 
Fixing the projection 
to $N$ is not very important but it is better to fix it.  

For example,  
first , on $\rho={\rho_{0}}$, 
namely on $\ol{U} \subset S^{2n+1}_{{\rho_{0}}}$, 
by the exponential map in $S^{2n+1}$ from $N$ we can define 
a projection to $N$. 
Then, for each ${\rho}'\in[{\rho_{0}}/2, {\rho_{0}}]$ 
the $CM\cap\{\ \rho={\rho}'\}$, 
applying the same construction on $S^{2n+1}_{{\rho}'}$, 
we obtain a smooth projection $CM \to N$. 
For the second step we can also rely on the gradient flow  
of the function $\rho$ on each $f^{-1}(w)$ ($w \in D^2_{\varepsilon_0}$) 
to extend the projection to $N$. 
(Also in the place of the gradient flow,  
the Liouville vector field on each $F_{w}=f^{-1}(w)$ is available.)  
This gives a smooth identification of $CM$ with 
 $[{\rho_{0}}/2, {\rho_{0}}]\times D^2_{\varepsilon_0}\times N$  
and 
we regard 
$(\rho, w, n)\in [{\rho_{0}}/2, {\rho_{0}}]\times D^2_{\varepsilon_0}\times N$ 
as the product coordinate of $CM$.   

Next take a smooth function $\psi(\rho)$ on $[{\rho_{0}}/2, {\rho_{0}}]$ 
satisfying 
\begin{eqnarray*}
&\bullet&\qquad 
\psi(\rho)\equiv 1, \qquad \rho \in [\frac{{\rho_{0}}}{2}, \frac{2{\rho_{0}}}{3}] , 
\qquad\qquad\qquad
\\
&\bullet& \quad 
0 \leq \psi(\rho)\leq 1, \quad \rho \in [\frac{2{\rho_{0}}}{3}, \frac{5{\rho_{0}}}{6}] , 
\\
&\bullet& \qquad 
\psi(\rho)\equiv 0, \qquad \rho \in [\frac{5{\rho_{0}}}{6}, {\rho_{0}}] . 
\end{eqnarray*}

With these preparations 
we show the following. 

\begin{prop}\label{modification}
{\rm  $\,$ 
For sufficiently small $\varepsilon_{0}$ and for any $w\in\C$ with 
$|w|= \varepsilon_{0}$, 
the identical embedding of $F_{w}$  
admits an isotopy through symplectic embeddings 
with respect to $\omega_0$ 
in such a way that on the part $\rho \leq {\rho_{0}}/2$ the isotopy is trivial, 
through the isotopy $\rho$ is preserved, and the part 
$F_{w}\cup\{\rho \geq \frac{5{\rho_{0}}}{6}\}$ 
is embedded onto  the corresponding part of 
$f^{-1}(0)$. 
Moreover these family of isotopies are chosen to be smooth on the parameter 
$\theta\in S^1$ and whole through the isotopies the projection 
to $N$ is stationary. }
\end{prop}

\noindent
{\it Proof}. \quad 
In the product coordinate, the hypersurface 
$F_{w}\cap \{{\rho_{0}}/2\leq \rho\leq {\rho_{0}}\}$ 
is considered to be a graph of the constant function 
$w$ from 
$[{\rho_{0}}/2, {\rho_{0}}] \times N$ to $D^2_{\varepsilon_0}$. 

As the modification of the identical embedding of $F_{w}$, 
take the function $\psi(\rho)w$ 
instead of the constant function $w$ and consider its graph. 
Let $s\in[0,1]$ be the parameter for the isotopy and 
consider the graphs of the smooth map 
$$
(\rho, n, w, s) \mapsto ((1-s)+  s\psi(\rho))w 
$$
from
$[{\rho_{0}}/2, {\rho_{0}}]\times N \times D^2_{\varepsilon_0} \times [0,1]$ 
to $D^2_{\varepsilon_0}$, 
where $(w, s)\in D^2_{\varepsilon_0} \times [0,1]$ are regarded 
as parameters so that the graphs are considered to be drawn in 
$[{\rho_{0}}/2, {\rho_{0}}]\times N \times D^2_{\varepsilon_0}$.  
For any $(w,s)$, on 
$[{\rho_{0}}/2, \frac{2{\rho_{0}}}{3}] \cup [\frac{5{\rho_{0}}}{6}, {\rho_{0}}]$ 
the graph is symplectic in $\C^{n+1}$. 
Also note that if $s=0$ or $w=0$, the 
graphs are symplectic.  
From the compactness  
we can see easily that there exists a positive 
$\varepsilon>0$ such that 
whenever $|w|\leq\varepsilon$ regardless to $s\in[0,1]$ 
the graphs is symplectic.  
To complete the proof take $w=\varepsilon_{0} e^{i\theta}$ 
for $\theta\in S^1$.  
\hspace*{\fill} $\square$
\\

Let $\Phi:M=\cup_{\theta\in S^1}F_{w} 
\cap D^{2(n+1)}_{\rho_{0}} \to \C^{n+1}$ 
be the map 
which is identity on $M\cap D^{2(n+1)}_{r/2}$ and 
is the above one with $w=\varepsilon_{0} e^{i\theta}$ and $s=1$.  
On each Milnor page, we take the restriction of 
$\Phi^*\omega_0$ 
and  
$\Phi^*\lambda_0$ 
which define a new exact symplectic structure on each Milnor page. 

Along the image of $\Phi$ take the symplectic orthonormal plane field 
to $\Phi(F_{w})$.  
It defines a symplectic parallel transport  
from the Milnor page 
$F_{w} 
\cap D^{2(n+1)}_{\rho_{0}}$, namely $\theta=0$,  
through the pages 
$f^{-1}(\varepsilon e^{i\theta}) 
\cap D^{2(n+1)}_{\rho_{0}}$ ($0\leq\theta\leq 2\pi$)  
to 
$F_{w} \cap D^{2(n+1)}_{\rho_{0}}$, $\theta=2\pi$,   
which gives rise to a symplectic monodromy.  
It is clearly stationary on 
$\frac{5r}{6}\leq {\rho} \leq {\rho_{0}}$  
so that the monodromy is identity around the boundary. 

Each page has the identical symplectically convex end, 
which also coincides with the end of 
$f^{-1}(0)\cap D^{2(n+1)}_{\rho_{0}}$. 
The boundary $N=\partial $, 
which is also called 
the binding of the open book 
or also the link of the singularity,  
 is 
naturally equipped with a standard contact structure 
$\xi_N=\mathrm{ker}\,\lambda_0\vert_N$.  
\\
\\

\begin{df}\label{df:ExactSymplOB}{\rm 
(Exact symplectic open book decomposition) $\,$ 
An open book decomposition of a $(2n+1)$-dimensional manifold 
${M}$ is an {\it exact symplectic} open book decomposition if 
each pages admits  exact symplectic structure with 
contact type boundary (or equivalently 
a globally convex symplectic structure (see \cite{EG})) 
and the monodromy 
preserves the symplectic structure and 
is identity in some neighborhood of the boundary.  
}
\end{df}
Then there exists a contact structure  $\xi$ on ${M}$, 
which is unique up to isotopy,  
such that the the binding is a contact submanifold which naturally 
coincides with one given as the boundary of the pages and   
$\xi$ is almost tangent to the pages.  
It is also said that the open book 
is supporting the contact structure $\xi$.  
For more details on these notions,  see \EG \cite{BCS}.

\begin{cor}\label{cor:ExactSymplOB}
{\rm  $\,$ 
With a Milnor open book of an isolated singularity, 
an exact symplectic open book decomposition is associated. 
}
\end{cor}

In our construction of symplectic foliations,  we need a 
less strict class of open books. 

\begin{df}\label{df:ConvexSymplOB}{\rm
(Convex symplectic open book decomposition)  $\,$ 
An open book decomposition of a $(2n+1)$-dimensional manifold 
${M}$ is a {\it convex symplectic} open book decomposition if 
each pages admits a symplectic structure with 
contact type boundary 
and the monodromy 
preserves the symplectic structure and 
is identity in some neighborhood of the boundary.  
Compared with the exact symplectic open books 
the condition is loosen from the exactness to not necessarily 
global convexity of the symplectic structures on pages. 

The binding $N$ of the open book which is regarded as 
the common boundary of all pages naturally admits a 
contact structure as the boundary of convex symplectic structure. 
The symplectic structure on a neighborhood of the boundary of a page 
is of the form of symplectization of the contact structure on $N$. 
}
\end{df} 
\subsection{Simple elliptic and cusp singularities}
\label{subsec:singularities}

\noindent
The singularities with which we are mainly concerned 
are the following isolated hypersurface singularities 
at the origin $0=(0,0,0)$
in three complex variables $(x,y,z)\in\C^{3}$. 
\begin{eqnarray*}
\widetilde{E}_ 6&\colon & x^3+y^3+z^3+kxyz, \; k^3\ne -27\\
\widetilde{E}_7 &\colon & x^2+y^4+z^4+kxyz, \; k^4\ne 64,\\
\widetilde{E}_8 &\colon & x^2+y^3+z^6+kxyz, \; k^6\ne 432. 
\end{eqnarray*}
These singularities are called {\it simple elliptic} singularities.  
Usually an simple elliptic singularity is defined as a 
normal surface singularity whose minimal resolution locus $E$ 
is a smooth elliptic curve. 
The self-intersection $E\cdot E$ 
is equal to $-3$, $-2$, and to $-1$ respectively for $\widetilde{E}_6$, 
 $\widetilde{E}_7$,  and for $\widetilde{E}_8$. 
The negative of  $E\cdot E$ is called the {\it index} of the 
simple elliptic singularity.     
It is known that the simple elliptic singularities which are 
realized as hypersurface singularities are analytically equivalent to 
one of the above three (Saito \cite{S}).  The restriction on $k$ is to 
ensure that the singularity is isolated, and $k$ can be zero.  
Even though the different $k$'s define 
analytically different singularities, they are, from smooth topological 
point of view, diffeomorphic to each other as   
the space of $k$ is connected.  
Therefore in this article we do not pay particular attention  to the value 
of $k$.

The following polynomial also defines an isolated singularity at the origin, 
which is called a {\it cusp} singularity:   
\begin{eqnarray*}
T_{pqr} \colon  x^p+y^q+z^r+kxyz, \; 
k \ne 0,\; \dfrac{1}{p}+\dfrac{1}{q}+\dfrac{1}{r}<1. 
\end{eqnarray*} 
In these cases, $k$ should just avoid $0$, then for different $k$'s 
they are linearly equivalent. 
For $k=0$ above polynomials define singularities of a different class.  
A cusp singularity is also usually defined as a normal surface singularity 
whose minimal resolution locus is a rational curve with a node or 
a certain cycle of smooth rational curves.  See, \EG \cite{H}. 
It is also know that a cusp singularity is analytically realizable 
as an isolated hypersurface singularity if and only if it is 
analytically equivalent to one of the above 
$T_{pqr}$-singularities for $1/p+1/q+1/r=1$ 
(Karras \cite{K}).     
\smallskip

The above three simple elliptic singularities can `almost' be considered 
as   $T_{pqr}$-singularities for $1/p+1/q+1/r=1$ 
with an attention to the value of $k$.  
The link of these singularities are known to be $T^{2}$-bundles 
over the circle with monodromy 
 $$A_{p,q,r}=\begin{pmatrix}r-1&-1\\1&0\end{pmatrix}
\begin{pmatrix}q-1&-1\\1&0\end{pmatrix}
\begin{pmatrix}p-1&-1\\1&0\end{pmatrix}.$$ 
See \cite{La}, \cite{N}, and \cite{Ka}.  
For simple elliptic singularities, \IE $(p,q,r)=(3,3,3)$, 
$(2,4,4)$, and $(2,3,6)$, 
$A_{p,q,r}$'s are conjugate to 
$\begin{pmatrix} 1&0 \\ 3&1\end{pmatrix}$, 
$\begin{pmatrix} 1&0 \\ 2&1\end{pmatrix}$ and 
$\begin{pmatrix} 1&0 \\ 1&1\end{pmatrix}$, respectively, 
and they are unipotent, so that the links are  nil 3-manifolds 
and also are considered to be the circle bundles over $T^{2}$ 
with euler numbers $-3$, $-2$, and $-1$ respectively (see \cite{Mi2}).  
For cusp singularities,  \IE  $1/p+1/q+1/r<1$,  
as $\mathrm{tr}\,A_{p,q,r} =2+pqr(1-(1/p+1/q+1/r))$, 
the monodromy is hyperbolic and the link is a solv 3-manifold.

\section
{Adapted foliations of open book decompositions}
\label{sec:Lawson}

In this section, 
 we review and discuss two criterions 
due to A. Mori (\cite{Mo}) for 
open books supporting contact structures
to be nicely foliated. 
The first one is on the binding so that 
the open book admits an adapted foliation 
and the second is on the contact structure of the 
binding so that the ambient contact structure 
converges to the foliation.  
Surprisingly these conditions give rise to 
leafwise symplectic foliations 
adapted to Milnor fibrations 
which we can associate starting from 
isolated singularities in complex three variables. 

\subsection
{Mori's criterions on open books
}
With an open book decomposition on a 3-manifold it is easy to 
construct an associated foliation of codimension one, 
which is called a {\it spinnable foliation} (cf. \cite{KMMM}).  

In general, foliating smoothly the exterior of the binding $N$ 
(which is denoted by $M$ in our context) is always easy.  
It suffices to make the pages winding around the boundary.  
However, foliating the tubular neighborhood of the binding 
is not trivial. 
Thurston's $h$-principle \cite{Th} 
tells us that for an open book of an odd dimensional closed manifold 
there always exists a foliation 
which has the boundary of the tubular neighborhood of the binding 
as a compact leaf, 
while from the $h$-principle we do not have enough control 
on the resultant foliation. 

\begin{thm}\label{Mori-1}
{\rm(Mori's criterion-1, \cite{Mo}, \cite{L}) $\,$ 
Let $M^{m}$ be a closed smooth manifold which is given an open book decomposition. 
If {\it the 
the binding $N^{m-2}$ of the open book admits a Riemannian 
foliation $\G$ of codimension one}, then 
its pull-back to $U\cong N\times {\mathrm{int}}D^2$ 
can be turbulized around the boundary $W=S^1 \times N$ 
so that  $\overline{U} \cong N\times D^2$ is smoothly 
foliated and the boundary $W$ is a compact leaf. 
As a conclusion  there exists a foliation 
$\F$ adapted to the open book decomposition.  
}\end{thm}

This is a slight generalization of Lawson's construction (\cite{L}) of 
foliation adapted to an open book, which just assumes that 
the binding fibers over the circle. Then we can foliate the tubular 
neighborhood of the binding by so called turbulization 
so as to have the boundary as a compact leaf. 
It is easy to choose the holonomy around this compact leaf 
so that after pasting the pages we obtain a smooth foliation. 
From the point of view of Mori's criterion, 
Lawson's construction adopts 
the bundle foliation which the fibration over the circle 
determines as the Riemannian foliation on the biding to start with.  
Such a foliation which is associated with 
such an open book is called a foliation of {\it Lawson type} or 
an {\it adapted} foliation.  
In this article, if we start from Mori's criterion-1with 
a Riemannian foliation, we call the resultant foliation 
a {\it generalized Lawson type} foliation.

As is well known, if the binding admits a Riemannian foliation of 
codimension one, a slight perturbation of the defining closed 1-form 
among the non-singular closed 1-forms so as to represent 
a rational first cohomology class, we obtain a Riemannian foliation 
determined by a fibration over the circle (Tischler \cite{Ti}). 
Therefore as the restriction of the topology of binding, 
Mori's criterion is equivalent to Lawson's one, however, 
classes of foliations on the binding and of resultant foliations 
are wider.  

Now we explain the turbulization.  First let us just assume 
the fibering of the binding $N$ over 
the circle $S^{1}=\{x\in\R/2\pi\Z\}$.   

Also fix a normal coordinate $w\in D^{2}_{\varepsilon_{0}}\subset\C$ 
of the fiber as the normal bundle to $N \subset U$ 
in the tubular neighborhood $\ol{U}$ and put $r=\vert w \vert$. 
Let $\psi(r)$ be a smooth function on  $[0, \varepsilon_{0}]$ 
satisfying 
\begin{eqnarray*}
&\bullet&\qquad \,
\psi(r)\equiv 1, \qquad r \in [0, \varepsilon_{0}/2] , 
\qquad\qquad\qquad
\\
&\bullet& \qquad 
\psi'(r) < 0, \qquad r \in (\varepsilon_{0}/2, \varepsilon_{0}),  
\\
&\bullet& \qquad 
\psi(\varepsilon_{0})=0, 
\\
&\bullet& \qquad 
\psi(r)\mathrm{\,\, is\,\, flat\,\, to\,\,} 0 \mathrm{\,\, at \,\,} \varepsilon_{0}.  
\end{eqnarray*}
Then take 
$$\alpha_{U}=\psi dx + (1-\psi)dr$$ 
as the defining $1$-form of the desired foliation on $\ol{U}$. 
\begin{rem}\label{rem:Reeb} {\rm $\,$ 
Consider the solid 3-torus $S^{1}\times D^{2}_{\varepsilon_{0}}$ 
with the coordinate $(x, w)$ and set $r=\vert w \vert$ as above.  
Then the $1$-form $\alpha$ as above defines a so called Reeb component, 
\IE a typical foliation on the solid torus. Two copies of Reeb components 
are glued together to form a Reeb foliation on $S^{3}$. 
The foliation on $\ol{U}$ constructed above 
can be regarded as the pull back of the Reeb component 
(or its open neighborhood) 
by the projection $p:  \ol{U} \to S^{1}$ ($p(n,w)=(x(p), w)$)
with the same monodromy 
and the fiber as those of $N\to S^{1}$.   

For a generalized Lawson type foliation, first take a defining closed 1-form 
$\alpha_{N}$ of the Riemanian foliation on $N$ and take a Tischler fibration 
to $S^{1}$. Then the place of $dx$ use $\alpha_{N}$ and take   
$$\alpha_{U}=\psi \alpha_{N} + (1-\psi)dr .$$  
Then it suffices to obtain a smooth foliation which is extendable 
even to the other side of the boundary as a trivial foliation 
defined by $dr$. 

The reason why the foliation with which we start should be 
Riemannian is as follows.  
If we trivially extend the defining $1$-form $\alpha_{N}$ 
of the foliation on $N$ to $\ol{U}=D^{2}\times N$ 
it also defines a foliation on $\partial U$.  
Due to Rosenberg and Thurston \cite{RT}, for the smooth turbulization 
which can be extended to the other side as a trivial foliation, 
it should be Riemannian.   
}
\end{rem}

The fibration of the binding over the circle works not only for foliating 
the tubular neighborhood of binding 
but also plays an important roll for existence of the leafwise 
symplectic structures.

\begin{thm}\label{Mori-2}
{\rm(Mori's criterion-2, \cite{Mo}) $\,$ 
Let $\xi$ be a contact structure on a closed manifold $M^{5}$ 
In the above criterion, \IE the existence of a Riemannian foliation 
on the binding $N^{3}$,  we also assume that 
{\it there exists a contact 1-form $\alpha$ defining $\xi$ which 
is adapted to the supporting open book, 
the contact form $\alpha\vert_N $ of 
the restricted contact structure $\xi|_{N}$ 
has the Reeb vector field is tangent to $\G$}. 
Then there exists a one-parameter family $\xi_t$ of 
contact structures, 
which is an isotopic deformation of $\xi$, 
converges to the Lawson type foliation $\F$ as hyperplane field.   
}\end{thm}

\begin{rem}\label{rem:criterion-2}{\rm  $\,$ 
This criterion applies to the Milnor open books of 
all simple elliptic and cusp singularities of hypersurfaces. 
However, the exact 2-forms $d\alpha$'s do not converge to 
the leafwise symplectic structures which are constructed 
in this article.  The compact leaf is an apparent obstruction. 

For the simple elliptic singularities, 
the natural contact structure obtained as the convex boundary 
of the Milnor page or the singular surface admits a Reeb vector 
field which is tangent to the circle fiber 
if the link is considered as a circle bundle over $T^{2}$.  
The circle fiber is contained in the $T^{2}$ fiber of the fibration 
to the circle. See \cite{Mi2}.  

In  the cusp case, 
a natural contact structure on the link which is a solv manifold 
admits a Reeb vector 
field tangent to one of the two eigen directions of the monodromy 
in each $T^{2}$ fiber.  See \cite{Ka} for this. 

Mori \cite{Mo} constructed a family of open books 
of $S^{4}\times S^{1}$ supporting contact structures  
whose bindings are nil 3-manifolds, 
where any of $-\ell$ for $\ell\geq 4$ are realized 
as the euler number of the link as circle bundle over $T^{2}$.    
}
\end{rem}

\subsection{Leafwise symplectic structure around the binding}
\label{subsec:AroundBinding}
In this subsection we explain the existence 
of leafwise symplectic structures on the tubular neighborhood of the 
binding. 

Surprisingly we will see the fibering of 
binding over the circle works further until 
the existence of the leafwise 
symplectic structures if we start from 
isolated singularities of complex three variables. 
\medskip

\noindent
In this subsection we assume 
\begin{itemize}
\item[(o)] the manifold $M^{2n+1}$ is closed and oriented,  
\item[(i)] $M$ admits an open book satisfying 
with a fixed fibration 
$p_{N}:N^{2n-1}\to S^{1}=\{x\in\R/2\pi\Z\}$ 
of the binding over the circle with the fiber $\Sigma^{2n-2}$ 
and the monodromy $\varphi_{\Sigma}\curvearrowright \Sigma$,  
\item[(ii)] the fiber $\Sigma$ is equipped with a symplectic structure 
$\omega_{\Sigma}$ and the monodromy  $\varphi_{\Sigma}$ preserves 
$\omega_{\Sigma}$,   
\end{itemize}
and let 
the fibration 
$p_{N}:N\to S^{1}$ 
trivially extend to the closed tubular neighborhood as 
$p:\ol{U}=N^{2n-1}\times D^{2}_{{\varepsilon_{0}}}\to S^{1}\times D^{2}_{{\varepsilon_{0}}}$. 

\begin{rem}\label{rem:binding}{\rm $\,$ 
To obtain a leafwise symplectic structure on a Lawson type foliation, 
we need that the compact leaf $\partial U=S^{1}\times N$ 
admits a symplectic structure.  
For this, the conditions (i) and (ii) are sufficient, because we can take 
$\omega_{\partial U}=dx \wedge d\theta + \tilde{\omega}_{\Sigma}$ 
where $\theta=\mathrm{arg}\,w$   
($w\in \partial D^{2}_{{\varepsilon_{0}}}
\subset\C$),  where    
$\tilde{\omega}_{\Sigma}$ is an extension of ${\omega}_{\Sigma}$ on a fiber 
to the total space of the fibration as a closed 2-form, 
which we have thanks to the condition (ii).  

For $n=2$, (i) implies (ii) 
because any orientation preserving mapping class of 
a closed oriented surface contains an area preserving representative.  
Moreover, (i) is also necessary for the symplectic structure on 
$\partial U$ due to Friedl-Vidussi \cite{FV}.    

Moreover, if we start from an isolated singularity of three complex 
variables ($n=2$), simple elliptic and cusp singularities are exactly 
those whose Milnor open book satisfy the condition (ii) (Neumann \cite{N}). 
For a simple elliptic singularity $\partial U$ is diffeomorphic to a 
so called Kodaira-Thurston nil-manifold which is also known as 
a Kodaira's primary surface.  
}
\end{rem}
\begin{thm}{\rm (Friedl-Vidussi, \cite{FV}) $\,$   
A closed oriented 4-manifold $W^{4} = S^{1}\times N^{3}$ 
admits a symplectic structure if and only if 
$N$ fibers over the circle. 
}
\end{thm}
In the case the symplectic structure on $W$ 
is not necessarily in the above form which we adopt in this article.

\begin{thm}{\rm (Neumann, \cite{N}) $\,$ 
The link of an isolated singularity of three complex variables 
is Seifert fibered or a graph manifold. If it fibers over the circle, 
it should be a nil 3-manifold or a solv 3-manifold, \IE a torus bundle 
over the circle, and the singularity is a simple elliptic one or a cusp.  
}
\end{thm}

\begin{constr}\label{constr:TubNbd}{\rm
(Leafwise symplectic structures on  $\ol{U}$, \cite{Mi2})  $\,$ 
We adopt the notations in the above remark. 
It is easy to take a leafwise symplectic structure $\omega_{R}$ 
on the Reeb component so that 
on the toral leaf $(\theta,\,x)\in S^{1}\times S^{1}$ it restricts to 
$d\theta \wedge dx$.  As the Lawson type foliation is the pull back 
of the Reeb component by the projection $p$, 
which is a trivial extension of $p:N \to S^{1}$ to the 
$D^{2}_{{\varepsilon_{0}}}$ factor, we also have an extension of 
$\tilde{\omega}_{\Sigma}$ which is denoted the same. 
Then, for any non-zero constants $a$ and $b$, 
$$
\omega_{\ol{U}} = a (d\theta \wedge dx) + b (\tilde{\omega}_{\Sigma})
$$ 
is a leafwise symplectic form which restricts to the boundary as 
desired above.  
}
\end{constr}
\section{Leafwise symplectic structures on adapted foliations}
\label{sec:MainTheorem}
We state the main results in this article and give their proofs .  
\begin{thm}\label{thm:MainTheorem}{\rm $\,$ 
Let $M$ be a closed oriented manifold of dimension $5$ 
equipped with a convex symplectic open book decomposition.  
Assume the following conditions. 
\begin{itemize}
\item[(1)] 
The binding $N$ admits 
a fibration $p_{N}:N\to S^{1}$ over the circle. 
Cnsequently, 
there is a closed 2-form $\tilde{\omega}_{\Sigma}$ on $N$ 
it is restricted to an area form on each fiber of $p_{N}$.   
\vspace{-5pt}
\item[(2)] 
The contact structure $\xi_{N}$ 
which the convex open book defines on $N$ 
admits a Reeb vector field which 
is tangent to the fibers of $p_{N}$.  
(Mori's criterion-2) 
\vspace{-5pt}
\item[(3)]
The cohomology class $[\tilde{\omega}_{\Sigma}]$ is 
in the image of 
the restriction map $H^{2}(F;\R) \to H^{2}(N;\R)$     
in the real cohomology 
of the Milnor fiber $F$ to its boundary $N\cong\partial F$.
This is equivalent to that 
the fundamental class $[\Sigma]$ of the fiber 
is non-trivial in  $H_{2}(F;\R)$.  
\vspace{-5pt}
\end{itemize}
Then the open book admits a leafwise symplectic foliation 
of Lawson type. 
}
\end{thm}
\begin{rem}\label{rem:condition(3)}{\rm $\,$ 
We can formulate the theorem in a slightly generalized way,  
by Mori's criterion-1 with a Riemannian foliation in the place of 
the foliation by the fibration $p_{N}$ in the conditions (1) and (2).   
In that case,  
the arguments in Subsection \ref{subsec:AroundBinding}
is fixed exactly as in Remark \ref{rem:Reeb}. 
For the arguments below we need to further assume that 
 $\tilde{\omega}_{\Sigma}\wedge\alpha_{\G}$  
 is a volume form of $N$ 
and to replace $dx$ with 
$\alpha_{\G}$,   
 where  $\alpha_{\G}$ is 
a closed 1-form on $N$ which defines the Riemannian foliation $\G$.

If we try to generalize these constructions 
to the Milnor open books of 
higher dimensional isolated singularities, 
for four or more complex variables, 
unfortunately the Milnor open book never satisfies not only 
(3) but also (1).  
That is because 
the Milnor fiber is 2-connected and the link is 1-connected.  
}
\end{rem}
\begin{cor}\label{cor:MainCorollary}{\rm $\,$ 
Any isolated hypersurface singularity of three complex variables 
which is simple elliptic or cusp 
has  
the Milnor open book of $S^{5}$ which admits 
leafwise symplectic foliation of Lawson type. 
}
\end{cor}
\noindent
{\bf Proof} of Corollay \ref{cor:MainCorollary}.  
For simple elliptic and cusp singularities of three complex variables, 
the condition (1) and the first half of (3) hold as explained in
1.3.   
The condition (2) is, as mentioned 
in Remark \ref{rem:criterion-2}
observed in \cite{Mi2} for the simple elliptic case 
and in \cite{Ka} for the cusp case.  
Therefore we need to check the condition (3).  
In \cite{Mi2} it is proved by  a Meyer-Vietoris argument 
in the simple elliptic case, 
while including that case it is easily shown  
from the 1-connectedness of the Milnor fiber $F$ 
is in this dimension (\cite{M}),   
a part of the cohomology long exact sequence for the pair 
$(F,N=\partial F)$ and the Poincare duality of $(F, N)$ as follows.    
$$
H^{2}(F;\R) \to H^{2}(N;\R) \to H^{3}(F, N;\R) 
\cong H_{1}(F;\R)=0. 
\vspace*{-5pt}
$$
\QED{\ref{cor:MainCorollary}}
\\
{\bf Proof} of Theorem \ref{thm:MainTheorem} goes along the 
three steps. The first two are preparatory, and in the third step 
we construct symplectic forms on the pages which are naturally 
glued smoothly to the leafwise symplectic structure 
on the closed tubular neighborhood $\ol{U}$. 
The construction is similar to that in \cite{Mi2}.  
\medskip
\\
\noindent
\underline{Step 1}.   
We put a better coordinates on the ends of pages 
by a standard procedure in symplectic geometry.  
Let 
$(F_{\theta}, \omega_{\theta}
$ 
($\theta\in S^{1}$) be the exact symplectic pages. 
The tubular neighborhoods of the boundaries of the pages are 
identified by the monodromy as symplectic manifolds 
so that we can assume they all are of the form
$(1-\varepsilon_{1}, 1]\times N, d(\varrho_{1}\alpha_{1}))$
of a part of symplectization of a contact structure
$\xi_{N}=\mathrm{ker}\,\alpha_{1}$ on the binding $N=\partial F_{\theta}$, 
where $\varepsilon_{1}>0$ is some constant and   
$\varrho_{1}$ is the coordinate of $(1-\varepsilon_{1}, 1]$.  
Once extend the end to $\varrho_{1}<\infty$ with 
the symplectic form $d(\varrho_{1}\alpha_{1})$ and take a 
contact form $\alpha_{N}$ for $\xi_{N}$ which satisfies the 
condition (4).   
Then, if necessary by taking  a constant multiplication, 
$\alpha_{N}$ can be regarded as the graph of the function 
$\varrho_{1}=\alpha_{N}/\alpha_{1}$ on $N$ with values in  
$[1, c]$ for some $c\geq 1$.   
Putting $\varrho\alpha_{N}=\varrho_{1}\alpha_{1}$  
and taking $\varrho=(\alpha_{1}/\alpha_{N})\varrho_{1}$ 
as the new coordinate function in the place of $\varrho_{1}$, 
the symplectic form is $d(\varrho\alpha_{N})$. 
Remark that in the case of Milnor open book $\varrho_{1} = \rho^{2}$.   
\medskip
\\
\noindent
\underline{Step 2}.   From the conditions (1) and  (2) the binding $N$ 
is equipped with a closed 2-form $\tilde{\omega}_{\Sigma}$ 
which restricts to a symplectic form on each fiber $\Sigma_{x}$. 
The preimage 
$r_{\theta}^{-1}(\tilde{\omega}_{\Sigma})\subset H^{2}(F;\R)$ 
is non-vacant affine space by the condition (3) where 
$r_{\theta}:  H^{2}(F;\R)\to H^{2}(N;\R)$ is the restriction map, 
even though the monodromy of the open book might act, 
it is not difficult to find a family of closed 2-forms $\kappa_{\theta}$ 
on each page $F_{\theta}$ satisfying the following conditions: 
\begin{itemize}
\item[(a)] On the end $\varrho\geq 1$  
$\kappa_{\theta}=\tilde{\omega}_{\Sigma}$,  
where by abuse of notation 
$\tilde{\omega}_{\Sigma}$ also denotes the pull back of 
$\tilde{\omega}_{\Sigma}$ by the projection 
$[1,\infty)\times N \to N$.
\item[(b)] On $\bigcup_{{\theta\in S^{1}}}\!F_{\theta} \subset M$, 
$\kappa\!=\!\{\kappa_{\theta}\}_{\theta\in S^{1}}$ is a smooth family of pagewise 
closed 2-forms.  
\end{itemize}
\medskip
\noindent
\underline{Step 3}.   
First remark that there is  a constant $b>0$ 
so that 
$||\omega_{\theta}^{2}||-b||\omega_{\theta}||\cdot ||\kappa_{\theta}||>0$ 
at any point of $F_{\theta}$.  
Consequently, 
$\omega_{1}=\omega_{\theta}+b\kappa_{\theta}$ 
is symplectic.  
For this  
we have to check $\omega_{1}^{2}>0$ 
on $\bigcup_{{\theta\in S^{1}}}\!F_{\theta} \setminus \{\varrho > 1\}$ 
and on  $\{\varrho \geq 1\}$.  
On  $\bigcup_{{\theta\in S^{1}}}\!F_{\theta} \setminus \{\varrho > 1\}$ 
the condition is open for $b$ and of course $b=0$ satisfies it. 
Therefore there exists $b_{0}>0$ such that $b$ with $0<b<b_{0}$ 
satisfies this open condition. 
On  $\{\varrho \geq 1\}$,  
for a product type Riemannian metric on $[1,\infty)\times N$ 
where we assume $||d\varrho||\equiv 1$,  
once we can verify 
$|| d\varrho\wedge\alpha_{N}\wedge d\alpha_{N}||
> 
b\cdot || d\varrho\wedge\alpha_{N} \wedge \kappa_{\theta}||$ 
on $\{\varrho=1\}$, 
the same inequality holds on   $\{\varrho \geq 1\}$.  
Remark that on  $\{\varrho \geq 1\}$ 
$\kappa_{\theta}={\tilde{\omega}}_{\Sigma}$.  
This inequality implies  $\omega_{1}^{2}>0$ on   $\{\varrho \geq 1\}$ 
because we have 
$(d\alpha_{N})^{2}=d\alpha_{N}\wedge{\tilde{\omega}}_{\Sigma}=0$.   
We prepare four positive constants 
\begin{eqnarray*}
{\mathbf a}
&=&
\min_{n \in N} || 
\alpha_{N}\wedge d\alpha_{N}||_{n},  
\qquad 
{\mathbf C}=\max_{n\in N}    
|| 
\alpha_{N} \wedge {\tilde{\omega}}_{\Sigma}||_{n}, 
\\
{\mathbf A}
&=&
\max_{n \in N} || 
\alpha_{N}\wedge d\alpha_{N}||_{n},  
\qquad 
{\mathbf m}
=
\min_{n \in N} ||dx \wedge \tilde{\omega}_{\Sigma} ||_{n}.   
\end{eqnarray*}
We take $b>0$ satisfying $b< {\mathbf a}/{\mathbf C}$ and  $b\leq b_{0}$.   
This is not only for that  $\omega_{1}$ is symplectic but also for the 
following arguments.

Now we deform the symplectic structures on the pages $F_{\theta}$, 
only on their ends which are identical independently of $\theta$, 
we omit $\theta$ from the description.   
On the end $[1,\infty)\times N$ the symplectic form is 
$\omega=d(\varrho\alpha_{N})$.  
As in \cite{Mi2} for some smooth functions 
$K(\varrho)$ and $L(\varrho)$, 
consider the following 2-form
$$
\omega_{E} =d(K(\rho)\wedge\alpha_{N})
+
b{\tilde{\omega}}_{\Sigma}
+
L(\varrho) d\varrho \wedge dx 
$$
 on $F\cap \{\varrho \geq 1\}$.   
Let us take $K(\rho)$ 
which satisfies the following conditions.  
\medskip
\\
\qquad (K-1) \quad  $K(\rho)=\rho$ \qquad \qquad  $1\leq \varrho \leq 3$,  
\\
\qquad (K-2) \quad   $0\leq K'(\rho)\leq 1$ \qquad  $3\leq \varrho <4$,  
\\
\qquad (K-3) \quad   $-1\leq K'(\rho)\leq 0$ \quad  $4< \varrho \leq 8$,   
\\
\qquad (K-4) \quad   $K(\rho)=0$ \qquad \qquad  $8\leq \varrho$.  
\medskip
\\
Then, we take $L(\varrho)$ 
which satisfies the following conditions.  
\medskip
\\
\qquad (L-1) \quad  $L(\rho)=0$ \qquad \qquad  $1\leq \varrho \leq 2$,  
\\
\qquad (L-2) \quad   $0<L'(\rho)\leq $  \qquad  $\,\,\,\,2\leq \varrho <3$,  
\\
\qquad (L-3)  \quad   $L(\rho)\equiv a$ \qquad \qquad  $3\leq \varrho$,  
\\
where $a$ is any constant satisfying 
$a>\dfrac{2{\mathbf A}+b{\mathbf C}}{b{\mathbf m}}$. 
\medskip
\\

On the end $[8, \infty)\times N$ the 2-form $\omega_{E}$ is 
a cylindrical symplectic form 
$ad\varrho \wedge dx + b{\tilde{\omega}}_{\Sigma}$ as desired, 
because in making the ends of pages winding around the boundary 
$\partial U=S^{1}\times N \ni (\theta, n)$ 
of the tubular neighborhood of the binding $N$, 
we take $\theta  = \varrho$  (mod $2\pi$).  

We have to check that $\omega_{E}$ is symplectic on 
$[1, 8]\times N$.   
The most important property of this construction is that 
at any point we have 
$$
d(K(\rho)\wedge\alpha_{N})
\wedge (d\varrho \wedge dx) =0\quad \cdots \quad(*).  
$$
This follows from the condition (2) (Mori's criterion-2) 
which implies $d\alpha_{N}$ is divisible by $dx$ 
as both of them annihilate the Reeb vector field of $\xi_{N}$, 
and just because  
$d(K(\rho)\alpha_{N})
=K'(\varrho)d\varrho\wedge\alpha_{N}+K(\rho)d\alpha_{N}$.   
Therefore  we have 
\begin{eqnarray*}
\omega_{E}^{2}
  &=&(d(K(\rho)\alpha_{N}))^{2} 
  + 2b(d(K(\rho)\alpha_{N}))\wedge{\tilde{\omega}}_{\Sigma}
  +2bL(\varrho)d\varrho \wedge dx\wedge{\tilde{\omega}}_{\Sigma}
\\
&=&2d\varrho\wedge\left(K'(\varrho)K(\varrho)\alpha_{N}
\wedge d\alpha_{N} 
+
bK'(\varrho)\alpha_{N}
\wedge{\tilde{\omega}}_{\Sigma}
+bL(\varrho) dx\wedge{\tilde{\omega}}_{\Sigma}
\right)
\end{eqnarray*} 
 and we also remark that $|K'(\varrho)K(\varrho)|<4$ holds. 
 
\begin{itemize}
\item On  $[1, 2]\times N$, $\omega_{P}=\omega_{1}$ is symplectic.    
\item On  $[2, 3]\times N$,  still 
$
d(K(\rho)\wedge\alpha_{N})
+
b{\tilde{\omega}}_{\Sigma}
=\omega_{1}$ holds and thus the sum of the first two terms 
is $\omega_{1}^{2}>0$ and the third term is also positive.  
\item On  $[3, 8]\times N$, the first two terms might be negative, 
while the third term is already big enough to conclude 
 $\omega_{E}^{2}>0$.    
\end{itemize}
On each page $F_{\theta}$, 
$\omega_{\theta}+b\kappa_{\theta}$ 
extends to the end to be $\omega_{E}$ as a symplectic form 
${\omega_{P}}_{\theta}$ 
and we obtain a leafwise symplectic structure 
${\omega_{P}}=\{{\omega_{P}}_{\theta}\}$  on 
$\bigcup_{\theta\in S^{1}}F_{\theta}$.  
Remark here that the subscript $\,{*}_{\theta}\,$ indicates 
the `page number'.  
This completes Step 3. 
\medskip

As is already remarked above, 
the end  $[8, \infty)\times N$ with the cylindrical 
(\IE $\varrho$-translation invariant) symplectic form 
$ad\varrho \wedge dx + b{\tilde{\omega}}_{\Sigma}$ 
is exactly the covering of 
$\partial U=S^{1}\times N$ with 
$ad\theta \wedge dx + b{\tilde{\omega}}_{\Sigma}$ 
by $\theta=\varrho$ (mod $2\pi$),  and naturally after the tubulization,  
by putting $\omega_{U}$ on $\ol{U}$  
and $\omega_{P}$ on $M\setminus \ol{U}$       
we naturally obtain a smooth leafwise symplectic structure 
after the tubulization (see \cite{Mi2} for more detail).   
\QED{\ref{thm:MainTheorem}}

\begin{rem}\label{rem:generalization}
{\rm  $\,$ 
For higher dimensional generalization, 
we need not only to assume (ii) 
in Subsection \ref{subsec:AroundBinding} 
but to put a stronger condition than (2) 
in order to achieve a good modification of 
the symplectic structures on the pages.   
We do not know good examples 
nor geometric way to state the extra condition 
which we need.  

If we take $\omega_{E}$ in the same way, 
\EG in the case $n=3$, in computing 
$\omega_{E}^{3}$,  we would have 
terms  
$bKL\, d\varrho\wedge dx 
\wedge d\alpha_{N}\wedge {\tilde{\omega}}_{\Sigma}$ 
or 
$bKK'\, d\varrho\wedge \alpha_{N} 
\wedge d\alpha_{N}\wedge {\tilde{\omega}}_{\Sigma}$ 
of which we do not know the positivity and thus 
they should be compared with 
$bL d\varrho\wedge dx 
\wedge {{\tilde{\omega}}_{\Sigma}}^{2}$.  
In the case of 
general dimensions there might appear terms 
\smallskip
\\
\qquad\qquad\qquad
$
b^{i}K^{n-i-1}L\,
d\varrho\wedge dx \wedge 
(d\alpha_{N})^{n-i-1}\wedge {{\tilde{\omega}}_{\Sigma}}^{i}
$
\vspace{-3pt}
\\
or 
\vspace{-3pt}
\\
\qquad\qquad\qquad
$b^{i}K^{n-i-1}K'\, d\varrho\wedge \alpha_{N} \wedge 
(d\alpha_{N})^{n-i-1}\wedge {{\tilde{\omega}}_{\Sigma}}^{i}
$
\smallskip
\\
for 
$
i=1,2,\cdots, n-2
$
and should be compared with 
$
b^{n-1}L \wedge d\varrho\wedge dx 
\wedge {{\tilde{\omega}}_{\Sigma}}^{n-i}
$. 
As we have to take $b$ small enough, 
we rather have to assume that the above terms vanish.  
Mori's criterion-2 kills the term of the above form of $i=0$ 
and in the case of $d=2$ we do not have other terms. 
So we would like to put the following condition
$$
 (d\alpha_{N})^{n-i-1} 
\wedge {{\tilde{\omega}}_{\Sigma}}^{i}=0
\quad \mathrm{on} \,\, N 
, \qquad i=1,2,\cdots, n-2. 
$$
Then all the arguments work exactly in the same way as in the case $d=2$ 
and we obtain a leafwise symplectic structure on the associated foliation  
with the open book.  }
\end{rem}

\section{A certain closed symplectic 4-manifold}\label{sec:ClosedSymplectic}
In this section, we present two constructions of closed 
symplectic manifolds as application of the construction 
in the previous section. 
\subsection{Gluing two Milnor fibers}\label{subsec:2fibers}

The first application 
of the end-cylindrical symplectic structure  
is a construction of a particular closed symplectic 4-manifold,   
which has similar properties as K3-surfaces,  
starting from a particular cusp singularity, $T_{2,3,7}$ or  $T_{4,4,4}$.    
\begin{prop}\label{prop:monodromy}{\rm $\,$ 
For the cusp singularity $T_{2,3,7}$ or  $T_{4,4,4}$, 
the monodromy of the link $N$ if the singularity 
is conjugate to its inverse as an element in $\mathit{SL}(2; \Z)$. 
}
\end{prop}
{\bf  Proof} of \ref{prop:monodromy}. The monodromy of the link of 
such a $T_{p,q,r}$-singularity is geievn as 
$A_{p,q,r}=\begin{pmatrix}r-1&-1\\1&0\end{pmatrix}
\begin{pmatrix}q-1&-1\\1&0\end{pmatrix}
\begin{pmatrix}p-1&-1\\1&0\end{pmatrix}$ 
as mentioned in Subsection \label{subsec:singularities}. 
Therefore we have 
$A_{2,3,7}=\begin{pmatrix}5&-11\\1&-2\end{pmatrix}$ 
and 
$A_{4,4,4}=\begin{pmatrix}21&-8\\8&-3\end{pmatrix}$ 
and they are conjugate to 
$\begin{pmatrix}2&1\\1&1\end{pmatrix}$ 
and 
$\begin{pmatrix}13&8\\8&5\end{pmatrix}$.   
As they are unimodular and symmetric, 
conjugate to their own inverses  
$\begin{pmatrix}1&-1\\-1&2\end{pmatrix}$ 
and 
$\begin{pmatrix}5&-8\\-8&13\end{pmatrix}$    
respectively.  
\QED{\ref{prop:monodromy}}
\\
\begin{constr}\label{constr:K3}{\rm\quad 
Let us take  $T_{2,3,7}$ and one of its Milnor fiber, \EG $F_{0}$ with 
the end-cylindrical symplectic structure $\omega_{P\,0}$ 
constructed in the proof of Theorem \ref{thm:MainTheorem} and 
consider the restriction $\omega(1)=\omega_{P\,0}|_{F(1)}$ 
to the truncated Milnor fiber $F(1)=F_{0}\cap\{\varrho\leq 10\}$. 
Let $(F(2), \omega(2))$ be the identical copy of $(F(1), \omega(1))$.  
We glue them along their boundaries. 
The boundary $N$ is the suspension 
of $T^{2}$ by the monodromy matrix $A_{2,3,7}\in \mathit{SL}(2; \Z)$, 
we can assume that 
the closed 2-form $\tilde{\omega}_{T^{2}}$ 
(here the fiber $\Sigma$ is $T^{2}$) is the `pul-back' 
by the suspension of the standard area form of $T^{2}$.  
The fact that the monodromy is conjugate to its 
inverse in  $\mathit{SL}(2; \Z)$ implies 
that even if we invert the orientation of the base circle, 
still it is isomorphic to the original  as oriented $T^{2}$-bundle. 
We identify the boundary $N(2)=\partial F(2)$ with the $N(1)=\partial F(1)$ 
by a map $\Phi : N(2)\to N(1)$ in this way, namely, 
the $T^{2}$-fiber orientation is preserved, but the orientation of the base circle 
is reversed as $x(1)\circ\Phi=-x(2) \in S^{1}=\R/2\pi\Z$. 
Here $x(i)$ denotes the $x$-coordinate which we have been fixed in previous   
sections on $N(i)=N$. Similar notations are used here for other coordinates.  
Therefore we have $\Phi^{*}dx(1)=-dx(2)$ while 
 $\Phi^{*}\tilde{\omega}_{T^{2}}=\tilde{\omega}_{T^{2}}$, and . 
 $\Phi : N(2)\to N(1)$ is orientation reversing. 
 Now in order to glue $F(1)$ and $F(2)$,  
around the boundaries take the identification of $\varrho$-coordinate as 
$\varrho(2)=20- \varrho(1)$.  
In this way  not only without changing the orientation as 4-manifolds 
but also without changing the symplectic structures around the boundaries,  
we can glue two copies $F(1)$ and $F(2)$ to obtain 
a closed symplectic manifold $W_{2,3,7}$, because  
we have $\Phi^{*}\tilde{\omega}_{T^{2}}=\tilde{\omega}_{T^{2}}$
and $d\varrho(1)\wedge dx(1)=(-d\varrho(2))\wedge(-dx(2))
=d\varrho(2)\wedge dx(2)$. 

Remark that this construction is 
not taking the `double'.  
Of course also to $T_{4,4,4}$ the same construction 
applies to obtain a closed symplectic 4-manifold $W_{4,4,4}$.  
}
\end{constr}

In both cases of $T_{2,3,7}$ and $T_{4,4,4}$ the Milnor fiber 
is simply connected and has Euler characteristic $12$.  
In fact it is not difficult to compute 
the Milnor number $\mu(T_{p,q,r})=p+q+r-1$ 
(see \cite{Ga}, \cite{KKMM}).  
Therefore the resultant 4-manifolds $W_{2,3,7}$ and $W_{4,4,4}$  
are simply connected and have the Euler characteristic 24.   
Also their Milnor lattice (the intersection form) is also computed 
(also see \cite{Ga}, \cite{KKMM}), and from these informations 
we see the integral cohomology ring of  $W_{2,3,7}$ 
and 
the rational cohomology ring of $W_{4,4,4}$  
coincide with those of a K3-surface respectively.  
In particular, we see $W_{2,3,7}$ is homeomorphic to a K3-surface. 
These cohomology computations and the arguments in 4-dimensional 
differential topology are informed to the author by M. Ue.  
Therefore it is strongly expected that they are diffeomorphic to 
a K3-surface.  
This strongly motivates the forthcoming paper \cite{KKMM}
which in fact achieved it.  

\subsection{Circular Example}\label{subsec:circular}
In this subsection, 
we slightly generalize the method in \underline{Step 3} of the proof 
of Theorem \ref{thm:MainTheorem}, 
namely a modification of symplectic structure on the end,  
and construct symplectic structures 
on $W=S^{1} \times N$.   
Here we present a construction for a Riemannian foliation on $N$ 
(see Remark \ref{rem:condition(3)}).  We start with the following 
objects.  
\medskip
\\
\quad $\bullet$  
$N$ : a $(2n-1)$-dimensional closed oriented manifold,  
\\
\quad $\bullet$  
$p_{N}:N\to S^{1}=\{x\in \R/2\pi\Z\}$ : a smooth fibration,   
\\
\quad $\bullet$  
$\tilde{\omega}_{\Sigma}$ : a closed 2-form on $N$ which 
restricts to a symplectic form on each
\\ 
\qquad\qquad \!fiber of $p_{N}$,     
\\
\quad $\bullet$  
$\xi_{N}$ : a positive contact structure on $N$,  
\\
\quad $\bullet$  
$\alpha_{N}$ : a contact 1-form defining $\xi_{N}$,
\\
\quad $\bullet$  
$\G$ : a Riemannian foliation of codimension one  on $N$,
\\
\quad $\bullet$  
$\alpha_{\G}$ : a non-singular 1-form which defines $\G$.   
\medskip
\\
We assume Mori's criterion 2 for the Reeb vector field 
of $\xi_{N}$ and $\G$.  We also assume the following condition.    
\medskip
\\
\qquad 
(*) : 
$(\tilde{\omega}_{\Sigma})^{n-1}\wedge\alpha_{\G}$ is 
a volume form on $N$. 
\medskip
\\
In the case that the Riemannian foliation $\G$ is chosen to be the 
foliation given by the fibration $p_{N}$, 
this condition is automatically satisfied  
(see Remark \ref{rem:condition(3)}). 
For the case $n\geq 3$ we further assume 
not only the condition (ii) in Subsection \ref{subsec:AroundBinding}
for a Tischler fibration $N\to S^{1}$ associated with the Riemannian 
foliation $\G$ but 
the codition 
in Remark \ref{rem:generalization} with $dx$.

In the case of $n=2$, for all the nil 3-manifolds 
$N=\mathit{Nil}^{3}(-l)$ with ($l\in \N$) 
and all the solv manifolds $N=\mathit{Solv}_{A}$ with hyperbolic monodoromy 
$A\in \mathit{SL}(2;\Z)$, namely, $\mathrm{tr}\,A\geq 3$, 
any of them admits a contact structure $\xi_{N}$ ]
and a contact 1-form $\alpha_{N}$ whose Reeb vector field 
is tangent to the fibers of the fibration $p_{N}:N\to S^{1}$.  
In these cases the condition is satisfied 
by taking the foliation by fibration 
as the Riemannian foliation $\G$.

Following the construction of symplectic form $\omega_{E}$ 
on the end of Milnor fiber in the proof of \ref{thm:MainTheorem}, 
let us consider the closed 2-form 
$$
\omega
=
d(K(\theta)\alpha_{N}) + \tilde{\omega}_{\Sigma} 
+L(\theta)d\theta\wedge \alpha_{\G}
$$
determined by smooth functions  $K(\theta)$ and $L(\theta)$ 
on $S^{1}$ and consider when 
$$
w^{n}=\left(K'(\theta)d\theta\wedge\alpha_{N} + K(\theta)d\alpha_{N}
 + \tilde{\omega}_{\Sigma} 
+L(\theta)d\theta\wedge \alpha_{\G}\right)^{n}
$$
gives an everywhere positive volume form.  
In the cohmology we have
$$[\omega]=[\tilde{\omega}_{\Sigma} ]
+
\left(\int_{S^{1}}L(\theta)\,d\theta\right)
[d\theta\wedge \alpha_{\G}]
\in H^{2}(W;\R).
$$

Because of Mori's criterion-2   
the term 
$d\theta \wedge \alpha_{\G}\wedge (d\alpha_{N})^{n-1}$ 
vanishes.  
Then as we have an extra condition that we assume for  the case $n\geq 3$,  
$w^{n}$ has only three terms as follows.  
\begin{eqnarray*}
& w^{n} &
\!\!=\,\, K'(\theta)K(\theta)^{n-1}d\theta \wedge \alpha_{N}
\wedge (d\alpha_{N})^{n-1} 
+ \,
L(\theta)d\theta\wedge\alpha_{\G}\wedge{\tilde{\omega}_{\Sigma}}^{n-1}
\\
&
&
\,\,+\,
 K'(\theta)d\theta \wedge \alpha_{N}
\wedge {\tilde{\omega}_{\Sigma}}^{n-1}.  
\end{eqnarray*}

We have a wide variety for the choices of  $K(\theta)$ and  $L(\theta)$ 
which make $\omega$  symplectic.  
For any choice of   $K(\theta)$, let us consider the 
ratio of the sun of the first and the third terms 
against 
$d\theta\wedge\alpha_{\G}\wedge{\tilde{\omega}_{\Sigma}}^{n-1}$ 
as follows
$$
\lambda(\theta)
=-\min_{n\in N} 
\frac{
K'(\theta)\alpha_{N}\wedge 
\left(K(\theta)^{n-1}
(d\alpha_{N})^{n-1} +
{\tilde{\omega}_{\Sigma}}^{n-1}\right)_{(\theta,n)}
}
{
{
\alpha_{\G}\wedge{\tilde{\omega}_{\Sigma}}^{n-1}
}_{(\theta, n)}
}
$$
where the ratio is taken as the that of $2n-1$ forms on $N$,   
abbreviating $d\theta$ in common. 
Then apparently $\omega$ is symplectic if and only if 
$L(\theta)$  satisfies 
$L(\theta)>\lambda(\theta)$ for any $\theta\in S^{1}$.   

Remark that the choice of  $K(\theta)$ 
and that of  $L(\theta)$ can be made quite independently.  
For example, such a choice of $K(\theta)$ and $L(\theta)$ 
is possible that $K(\theta)\equiv 0$ on some intervals while 
 $L(\theta)\equiv 0$ on some other intervals.

\begin{thm}\label{thm:circular}
{\rm $\,$ 
For a closed 3-dimensional solv 
or nil manifold $N$, $W=S^{1}\times N$
admits a symplectic form 
which looks like a co-symplectic form
(see \ref{subsec:b-symplectic}) 
on some intervals in $S^{1}$ and 
like a slight deformation of 
a symplectization of the canonical contact structure 
on $N$ on some other intervals. 
}
\end{thm}

\section{Other applications}\label{sec:other}

\subsection{$\boldsymbol{b}$-Symplectic structures}\label{subsec:b-symplectic}
As another application of the construction 
of end-cylindrical symplectic structure, 
we remark in this section that certain $b$-symplectic manifolds 
are easily constructed.   

A {\it b-symplectic} structure is a Poisson structure (Poisson bi-vector field) 
$\Pi\in C^\infty(W; \bigwedge^2TW)$ 
on a $2n$-dimensional manifold $W$ such that the graph of 
$\wedge^n\Pi$ is transverse to the zero section 
in $C^\infty(W; \bigwedge^nTW)$. 
Therefore the singular locus $Z=\{w\in W\,\vert\, 
\wedge^n\Pi(w)=0\}$ is a smooth hypersurface 
of $W$. 
This notion was introduced by Melrose in \cite{Me} 
and recently is drawing more attensions from Poisson geometry. 
For more details, see \EG \cite{GMP}.  

For a $b$-symplectic manifold,  
a Darboux theorem holds around the singular locus $Z$ 
in the following a way.   
For a point in $Z$ 
there exists a coordinate neighborhood 
$(x_{1}, y_{1}, \ldots, x_{n}, y_{n})$ so that 
the Poisson bi-vector $\Pi$ is expressed as
$$
\Pi=x_{1}\frac{\partial}{\partial x_{1}}\wedge
\frac{\partial}{\partial y_{1}}
+
\sum_{i=2}^{n}
\frac{\partial}{\partial x_{i}}\wedge
\frac{\partial}{\partial y_{i}}
$$
where $Z=\{x_{1}=0\}$.  
This is equivalent to say that the corresponding symplectic form 
$\omega$ has the expression   
$$
\omega=\frac{1}{x_{1}}dx_{1}\wedge dy_{1}
+
\sum_{i=2}^{n}
dx_{i}\wedge dy_{i}.
$$
If we allow the degeneration (or the divergence) along the 
codimension one singular locus $Z$ to be the $\ell$-th order 
$\ell\in\N$ in transverse direction, 
it is called a $b^{\ell}$-{\it symplectic structure} 
and has a local expression replacing the coefficient 
$x_{1}^{\pm 1}$ above with $x_{1}^{\pm\ell}$ 
($\pm$ correspond to Poisson/symplectic formulations).  
\medskip

Let $W$ be the double of 
the Milnor page of one of simple elliptic or cusp singularities. 
Namely, take two copies of a compact Milnor page $F_\theta$  
with boundary and paste them along the boundaries with identity 
to obtain 
$W=-F_\theta \cup_\partial F_\theta$.  

\begin{thm}\label{thm:b-symplectic}{\rm   $\,$ 1) $\,$
$W$ admits a $b$-symplectic structure whose  
singular locus is 
$\partial F_{\theta}=N$. 
\\
2) $\,$
Exactly the same holds for  $b^{\ell}$-symplectic structures 
for odd $\ell\in\N$. 
\\
3) $\,$ For any even $\ell\in\N$,  
$W_{2,3,7}$ and $W_{4,4,4}$ in Subsection \ref{subsec:2fibers} 
admit $b^{\ell}$-symplectic structures  whose  
singular loci are    
$\partial F(i)=N$. 
}
\end{thm}

Now our purpose in this section is a consequence of 
the following rather general remark . 

\begin{prop}\label{prop:b-symplectic}{\rm  $\,$ 
Let  $(W_{\pm}, \Omega_{\pm})$  be  
a pair of $2n$-dimensional symplectic manifolds 
with boundary which have 
the identical collar neighborhoods of the boundaries 
as follows.  \vspace{-5pt}
\begin{enumerate}
\item[(1)] $\partial W_+$ and $\partial W_-$ admit collar neighborhoods 
both of which are symplectomorphic to $([0,2)\times N, \Omega)$.   
\vspace{-5pt}
\item[(2)] There exist a non-singular 
closed $1$-form $\alpha_{N}$ 
and a closed $2$-form $\Omega_N$ 
on $N$ satisfying \vspace{-5pt}
$$
\Omega=dt\wedge{P_{N}}^*\alpha_{N}+ {P_{N}}^*\Omega_{N}
\vspace{-5pt}
$$
where ${P_{N}}$  denotes the projection  $[0,2)\times N \to N$ 
and $t$ denotes the coordinate on $[0,2)$.  \vspace{-5pt}
\end{enumerate}
Then, the manifold $W=W_-\cup_{\partial W_-=\partial W_+}W_+$ 
obtained by identifying the boundaries admits a 
$b$-symplectic structure whose singular locus is 
$\partial W_-=\partial W_+$.   
}
\end{prop}

\begin{rem}{\rm $\,$ 1) $\,$ 
The condition around the boundary is also said that the boundary 
admits a {\it cosymplectic structure}. 
More precisely, if a compact boundary $N$ 
(or a compact closed hypersurface) 
of a symplectic manifold $(W, \Omega)$ admits a 
closed one form $\alpha_{N}$ with 
$\alpha_{N}\wedge {\Omega_N}^{n-1}$ a volume form of $N$, 
then $N$ is called of {\it coymplectic type}, 
and the pair $(\alpha_{N}, \Omega_N)$ is called a cosymplectic 
structure. 
From the normal form theorem due to Gotay \cite{Go}, 
$N$ has a collar neighborhood where $\Omega$ is described 
exactly as in the above Proposition. 
For our application we can directly arrange the situation in 
a stronger form as in the proposition. 
\smallskip
\\
2) $\,$ 
For $\ell\geq2$, the claim on  $b^{\ell}$-symplectic structures is almost direct. 
Even for $\ell=0$ the claim is nothing but Construction \ref{constr:K3}. 
\smallskip
\\
3) $\,$ 
The assumption on such a closed 1-form $\alpha_{N}$ 
is nothing but Mori's criterion-1.  
}
\end{rem}
\noindent
{\bf Proof} of Proposition \ref{prop:b-symplectic}. \quad 
Let $t_{\pm}\in [0,2)$ 
denote the each coordinate for $W_{\pm}$ corresponding to $t$.   
For the resultant manifold $W$, let $\tau\in (-2,2)$ be the corresponding 
coordinate and we construct $W$ by gluing  $W_{\pm}$ by defining 
$\tau=-t_{-}$ on $(-2,0]$ 
and 
$\tau=t_{+}$ on $(-2,0]$.   
The two symplectic forms 
$\Omega_{\pm}=\Omega\,\,\mathrm{on}\,\,W_{\pm}$ 
do not coincide at $\tau=0$ 
because they are indicated as 
$\Omega_{\pm}= \pm dt\wedge{P_{N}}^*\alpha_{N}+ {P_{N}}^*\Omega_{N}$. 
Remark that the orientations of $W_{\pm}$ are opposite to each other 
in $W$.  
We do not need to orient $W$ but if we define it as that of $W_{+}$, 
the original orientation of $W_{-}$ is negative. 

In order to obtain a $b$-symplectic structure whose singular locus 
is $N\times\{0\}$,  
Take a smooth positive function $p(\tau)$ 
on $(-2,2)$ satisfying
\medskip
\\ 
\qquad \qquad(i)  \,\,\,\,
$p(\tau)=1/t$\qquad$\tau\in (-1/2,0)\cup (0,1/2)$, 
\\ 
\qquad \qquad(ii)  \,\,
$p(\tau)\equiv -1$ \qquad$\tau\in (-2, -3/2)$, 
\\ 
\qquad \qquad(iii)  \,
$p(\tau)\equiv 1$ \quad\qquad$\tau\in (3/2, 2)$, 
\\ 
\qquad \qquad(iv)  \,
$p'(\tau)< -0$ \qquad$\tau\in (-3/2, -1/2)\cup(1/2,3/2)$.   
\medskip
\\
Then replace $\Omega$ with 
$\Omega_b=p(t)dt\wedge{P_{N}}^*(\alpha_{N}) + {P_{N}}^*\Omega_N$ 
on $(-2,2)$ and

the collar neighborhoods and put them on $W_i$'s. 
They naturally give rise to a $b$-symplectic structure on $W$.  
\hspace*{\fill}\QED{\ref{prop:b-symplectic}}
\\
\subsection{Leafwise symplectic foliations}\label{subsec:LSF}
Taking the double 
is useful in constructing regular leafwise symplectic foliations 
of codimension one.  Because it concerns only the cylindrical 
end of the symplectic manifold, let us start only with 
the cylindrical end of the setting of Proposition \ref{prop:b-symplectic} 
and follow some notations there.  
\begin{constr}\label{constr:double}{\rm  $\,$ 
Let 
$
W_{a}=(a, \infty) \times N$ ($a$ has no particular meaning)
be a 
$2n$-dimensional symplectic manifold 
with a symplectic form 
$$
\Omega=dt\wedge{P_{N}}^*\alpha_{N} + {P_{N}}^*\Omega_N
$$
as in Proposition \ref{prop:b-symplectic} (2).   
We construct a symplectic foliation of codimension one 
on $M'=
(-2,2)
\times N \times S^{1}$. 
Take a smooth function $\varphi(\tau)$ on $(-2,2)
$ 
satisfying 
\begin{eqnarray*}
& \mathrm{(i)}\,\,\varphi(\tau)\equiv 0 \,\,\,\, (\tau\leq-1
), \qquad\quad\quad
\mathrm{(ii)}\,\, \varphi(\tau)\equiv \pi \,\,\,\, (\tau\geq 1
)\qquad\quad
\\
& \mathrm{(iii)}\,\, \varphi(0)=\frac{\pi}{2} \, , 
\qquad\qquad\qquad\qquad
\mathrm{(iv)}\,\, 
\varphi'(\tau)\geq0\,\,\,\,(-1<\tau<1), 
\\
& \mathrm{(v)}\,\, 
\varphi(\tau)\equiv \frac{\pi}{4} \,\,\,\,(-\frac{2}{3}<\tau<-\frac{1}{3}), 
\quad
\mathrm{(vi)}\,\, 
\varphi(\tau)\equiv \frac{3\pi}{4} \,\,\,\,(\frac{1}{3}<\tau<\frac{2}{3}), 
\end{eqnarray*}
and consider the hyperplane field $\tau \F'=\mathrm{ker}\,\,\alpha'$ 
defined by the 1-form 
$\alpha'=\cos \varphi(\tau) d\theta -\sin \varphi(\tau) dt$ 
where $\tau$ and $\theta$ are the coordinates of 
$(-2, 2)$ and $S^{1}$ respectively and extend them 
to.  
Apparently the hyperplane field is integrable and 
defines a foliation $\F'$ of codimension one which is tangent to 
`$\mathrm{arg}=\varphi$' in $\tau\theta$-plane.  
$M'$ is covered  
by three open sets 
\medskip
\\
\qquad\qquad
$M_{-}=\{-2<\tau<-1/3\}$,  
\\
\qquad\qquad
$M_{0}\,=\{-2/3<\tau<2/3\}$, and 
\\
\qquad\qquad
$M_{+}=\{1/3<\tau<2\}$.  
\medskip
\\
We define  
the leafwise symplectic form $\Omega'$ on $(M', \F')$ 
by restricting the 2-forms 
\smallskip
\\
\qquad\qquad
$
\Omega_{-}=\,\,\,\, d\tau\wedge{P_{N}}^*\alpha_{N}
+
{P_{N}}^*\Omega_N    
$  
\qquad 
on $M_{-}$, 
\\
\qquad\qquad
$
\Omega_{0}\,=\,\,\,\,\, d\theta\wedge{P_{N}}^*\alpha_{N}
+
{P_{N}}^*\Omega_N    
$  
\qquad
on $M_{-}$, and 
\\
\qquad\qquad
$
\Omega_{+}=-d\tau\wedge{P_{N}}^*\alpha_{N}
+
{P_{N}}^*\Omega_N    
$  
\qquad 
on $M_{+}$  
\medskip
\\
to the leaves.  
Thanks to the conditions (v) and (vi), on 
$M_{-}\cap M_{0}$ and on $M_{0}\cap M_{+}$ 
the two restrictions coincide to each other.  
Remark that the restriction of $(M',\F', \Omega')$ to 
$M'\cap\{\tau<0\}$ and $M'\cap\{\tau>0\}$ are both 
isomorphic to 
$(W_{a}\times S^{1}, \,\,
\F=\{W_{a}\times \{\theta\in S^{1}\}\}, \,\,
\Omega
=dt\wedge{P_{N}}^*\alpha_{N} + {P_{N}}^*\Omega_N)
$  
as leafwise symplectic foliations.  
}\end{constr}

As we have seen, we can modify the symplectic structures 
of the Milnor fibers of the simple elliptic or cusp 
singularities of hypersurface type so as to have such ends.  

\begin{thm}\label{thm:double}{\rm  $\,$ 
Let $F_{\theta}$ be the Milnor fiber of one of 
$T_{p,q,r}$ singularities for $1/p+1/q+1/r < 1$ 
or $\widetilde{E}_k$ for $k=6,7,8$,  
and $W=-F_{\theta}\cup_{\partial F_{\theta}}\!F_{\theta}$ 
be its double.   
Then, $W\times S^{1}$ admits a leafwise symplectic foliation 
of codimension one.
 }
 \end{thm}

\begin{flushright}
Yoshihiko MITSUMATSU 
\\
\mbox{}
\\
{\small Department of Mathematics, 
\\
Chuo University
\\
1-13-27 Kasuga Bunkyo-ku, 
\\
Tokyo, 112-8551, 
\\
Japan
\\
\mbox{}
\\
yoshi@math.chuo-u.ac.jp
\\
}
\end{flushright}

\end{document}